\documentclass[12pt]{article}
\usepackage{amsmath}
\usepackage{amssymb}
\usepackage{amsopn}
\usepackage{latexsym}
\usepackage{amsfonts}
\title{\textbf{Renormalised solutions in thermo-visco-plasticity for a Norton-Hoff type model. Part I: the truncated case}}
\author{\textbf{Krzysztof Che{\l}mi\'{n}ski}\\[0.5ex]
\textbf{\footnotesize{Faculty of Mathematics and
 Information Science, Warsaw University of Technology,}}\\[-1ex]
\textbf{\footnotesize{ul. Koszykowa 75, 00-662 Warsaw, Poland}}\\[-1ex]
\textbf{\footnotesize{E-Mail: kchelmin@mini.pw.edu.pl}}\\[1ex]
\textbf{Sebastian Owczarek}\\[0.5ex] 
\textbf{\footnotesize{Faculty of Mathematics and
 Information Science, Warsaw University of Technology,}}\\[-1ex]
\textbf{\footnotesize{ul. Koszykowa 75, 00-662 Warsaw, Poland}}\\[-1ex]
\textbf{\footnotesize{E-Mail: s.owczarek@mini.pw.edu.pl}}}
\date{}
\newtheorem{tw}{Theorem}[section]
\newtheorem{lem}[tw]{Lemma}
\newtheorem{de}[tw]{Definition}

\DeclareMathOperator{\dev}{dev}

\begin{document}
\maketitle
\begin{abstract}
We prove existence of global in time strong solutions to the truncated thermo-visco-plasticity with an inelastic constitutive function of Norton-Hoff type. This result is a starting point to obtain renoramlised solutions for the considered model without truncations. The method of our proof is based on Yosida approximation of the maximal monotone term and a passage to the limit.
\end{abstract}
\newcommand{\bl}{\backslash}
\newcommand{\nn}{\nonumber}
\newcommand{\KK}{\sigma_{\rm y}}
\newcommand{\ve}{\varepsilon}
\newcommand{\R}{{\mathbb R}}
\newcommand{\D}{{\mathbb C}}
\newcommand{\E}{{\cal E}}
\newcommand{\K}{{\cal K}}
\newcommand{\di}{{\mathrm d}}
\renewcommand{\S}{{\cal S}^3}
\renewcommand{\SS}{{\cal S}^3_{\mathrm{dev}}}
\newcommand{\id}{ {1\!\!\!\:1 } }
\section{Formulation of the problem and main result}
Our study is directed to mathematical analysis of thermo-visco-plasticity (for derivation we refer to \cite{9,25} and \cite{11}). This means to problems from the theory of inelastic deformations in which the temperature affects the visco-plastic response of the considered material. Let us assume that the elastic constitutive stress-strain relation has the form
\renewcommand{\theequation}{\thesection.\arabic{equation}}
\setcounter{equation}{0}%
\begin{eqnarray}
\label{eq:1.7}
\sigma&=&\D(\ve(u)-\ve^p)-f(\theta)\id\,,
\end{eqnarray}
where $\sigma$ is the Cauchy stress tensor, $\ve(u)=\frac{1}{2}(\nabla u +\nabla^T u)$ is the linearized strain tensor, $u$ is the displacement vector, $\ve^p$ is the inelastic strain tensor, $\theta$ is the temperature, f is a given function depending on the considered material and $\D$ is the elasticity tensor which we assume to be symmetric and positive definite on the space of symmetric matrices. Notice that the thermal part of the stress $-f(\theta)\id$ is not linearized in the neighborhood of the reference temperature (compare \cite{1}-\cite{3} and \cite{12,13}). Our motivation for the form of the elastic constitutive relation (\ref{eq:1.7}) follows the results of \cite{6} and \cite{7}.\\
Moreover, we assume that the density of the internal energy $e$ has the simple form\\ $e=c\,\theta+\D^{-1}T\cdot T$, where $T=\sigma+f(\theta)\id$. In this case as a consequence of the first principle of thermodynamics we obtain the following form of the heat equation 
\begin{eqnarray}
\label{eq:1.8}
c\rho\,\theta_t-\varkappa\Delta\theta+f(\theta)\mathrm{div}\,u_t&=&T\cdot\ve^p_t\,,
\end{eqnarray}
where $\rho>0$ is the mass density, $\varkappa>0$ is the material's conductivity. It is also assumed that the evolution of the inelastic strain tensor $\ve^p$ is given in the form
\begin{eqnarray}
\label{eq:1.9}
\ve^p_t=G(T),
\end{eqnarray}
where $G$ is a given maximal monotone vector field with $G(0)=0$ (inelastic constitutive equation with only one internal variable $\ve^p$, see \cite{18}). In the main part of the article we will specify $G$ choosing it in a form of the Norton-Hoff model (similar nonlinear flow rule was considered in \cite{14}).\\
If we consider equation (\ref{eq:1.8}) with the homogeneous Neumann boundary condition and the homogeneous balance of forces $\mathrm{div}\, \sigma=0$ also with homogeneous boundary conditions then we can conclude that the considering problem possesses a natural semi-invariant function namely the total energy
\begin{eqnarray}
\label{eq:1.10}
\E(t)=\int_{\Omega}c\rho\,\theta\,\di x+\frac{1}{2}\int_{\Omega}\D^{-1}T\cdot T\,\di x\leq\E(0),
\end{eqnarray}
where $\Omega\subset\R^3$ is a bounded domain describing our considered body with boundary of class $C^2$. In this simple calculation we have used the dissipation inequality, which yields that $T\cdot\ve^p_t\geq 0$ in the whole deformation process. From the observation (\ref{eq:1.10}) we see that the temperature is controlled in the space $L^1(\Omega)$ only (with additional information that $\theta\geq 0$). Moreover, in general the term $T\cdot\ve^p_t$ belongs to $L^1(\Omega)$ only. From these reasons we are going to prove existence of solution in renormalised sense -see for example \cite{31}-\cite{5} and \cite{7,20,19,16,32}. In the literature there are only some articles with mathematical study of special thermo-visco-plastic models with various modifications (see for example \cite{1}-\cite{17} and\cite{8,10}). See also the articles \cite{21} and \cite{22}, where a poroplasticity models
are investigated which have a similar structure to the linear thermo-plasticity.\\
In this article and in the following work \cite{15}, we are going to study the therm-visco-plasticity model of the Norton-Hoff type with a damping term, which we interpret as external forces acting on the material and depending on the deformation velocity. Thus the system of equations, which we study in this article is in the form
\begin{eqnarray}
\label{eq:1.11}
\mathrm{div}_x\, \sigma&=&-F-\mathrm{div}\,\D (\ve(u_t))\,,\nn\\
\sigma&=&\D(\ve(u)-\ve^p)-f(\theta)\id\,,\nn\\[1ex]
\ve^{p}_{t}&=&|\dev(T)|^{p-1}\dev(T)\,,\\[1ex]
T&=&\D(\ve(u)-\ve^p)\,,\nn\\[1ex]
\theta_t-\Delta\theta+f(\theta)\mathrm{div}\,u_t&=& |\dev(T)|^{p+1}\,,\nn
\end{eqnarray}
where $p>1$ is a given real number and $\dev\,(T)=T-\frac{1}{3}\,\mathrm{tr}\,(T)\cdot\id$ is a deviatoric part of the symmetric matrix $T$. The function $F:\Omega\times[0,T]\rightarrow \R^3$ describes the density of the applied body forces. The main idea in the existence theory of renormalised solutions is to study the so called truncated problem and next to prove that the sequence of obtained solutions converges to a renormalised solution.\\
In this article we study the truncated model only, hence for fixed $T>0$ and $\epsilon>0$ we have to find the displacement field $u:\Omega\times[0,T]\rightarrow \R^3$, the temperature of the material $\theta:\Omega\times [0,T]\rightarrow \R$ and the visco-plastic strain tensor $\ve^p:\Omega\times [0,T]\rightarrow\SS$ ($\SS$ denotes the set of symmetric $3\times 3$-matrices with vanishing trace) satisfying the following system of equations
\begin{eqnarray}
\label{eq:1.1}
\mathrm{div}_x\, \sigma&=&-F-\mathrm{div}\,\D (\ve(u_t))\,,\nn\\
\sigma&=&\D(\ve(u)-\ve^p)-f\big(T_{\frac{1}{\epsilon}}(\theta))\id\,,\nn\\[1ex]
\ve^{p}_{t}&=&|\dev(T)|^{p-1}\dev(T)\,,\\[1ex]
T&=&\D(\ve(u)-\ve^p)\,,\nn\\[1ex]
\theta_t-\Delta\theta+f\big(T_{\frac{1}{\epsilon}}(\theta)\big)\mathrm{div}\,u_t&=&T_{\frac{1}{\epsilon}}(|\dev(T)|^{p+1})\,,\nn
\end{eqnarray}
where the function $T_{\frac{1}{\epsilon}}(\cdot)$ is the truncation at height $\frac{1}{\epsilon}>0$ i.e. $T_{\frac{1}{\epsilon}}(r)=\min(\frac{1}{\epsilon},\max(r,-\frac{1}{\epsilon}))$. In this paper we assume that the function $f\in C^1(\R;\R)$. The system (\ref{eq:1.1}) is considered with the Dirichlet boundary condition for the displacement and the Neumann boundary condition for the temperature
\begin{eqnarray}
\label{eq:1.2}
u(x,t)&=&g_D(x,t)\quad  \textrm{ for}\quad x\in \partial \Omega \quad\textrm{and}\quad  t\geq 0\,,\nn\\[1ex]
\frac{\partial\,\theta}{\partial\,n}(x,t)&=&g_{\theta}(x,t)\quad\, \textrm{ for}\quad x\in \partial \Omega \quad\textrm{and}\quad  t\geq 0\,.
\end{eqnarray}
Finally, we consider the system (\ref{eq:1.1}) with the following initial conditions
\begin{eqnarray}
\label{eq:1.3}
u(x,0)=u_0(x),\quad \ve^p(x,0)=\ve^{p,0}(x),\quad \theta(x,0)=\theta_0(x).
\end{eqnarray}
Suppose that our data have the following regularity
\begin{equation}
\label{eq:1.4}
F\in H^{1}(0,T;L^2(\Omega;\R^3))\,,\quad g_D\in W^{2,\infty}(0,T;H^{\frac{1}{2}}(\partial\Omega;\R^3))\textrm{,}\\[1ex]
\end{equation}
\begin{equation}
\label{eq:1.5}
g_{\theta}\in H^{1}(0,T;L^2(\partial\Omega;\R))\textrm{,}\\[1ex]
\end{equation}
\begin{equation}
\label{eq:1.6}
u_0\in H^2(\Omega;\R^3) \textrm{,}\quad \ve^{p,0}\in L^2(\Omega;\SS)\textrm{,}\quad \theta_0\in H^1(\Omega;\R)\,.
\end{equation}
Next, we define a notion of a solution for the system (\ref{eq:1.1}).
\begin{de}
\label{de:1.1}
Suppose that the given data satisfy (\ref{eq:1.4})-(\ref{eq:1.6}). We say that a vector\\ $u\in H^{1}(0,T;H^1(\Omega;\R^3))$, the function $\theta\in L^{\infty}(0,T;H^1(\Omega;\R))$ such that $\theta_t\in L^{2}(0,T;L^2(\Omega;\R))$ and the inelastic deformation tensor $\ve^p\in H^{1}(0,T;L^2(\Omega;\SS))$ are solutions of the problem (\ref{eq:1.1})-(\ref{eq:1.3}) if the equations $(\ref{eq:1.1})_1$ - $(\ref{eq:1.1})_4$ are satisfied for almost all $(x,t)\in\Omega\times (0,T)$ and the heat equation is satisfied in the following sense 
\begin{eqnarray*}
\int_{\Omega}\theta_t\,v\,\di x+\int_{\Omega}\nabla\theta\nabla v\,\di x&+& \int_{\Omega} f\big(T_{\frac{1}{\epsilon}}(\theta)\big) \mathrm{div}\,u_t\,v\,\di x =\int_{\Omega}T_{\frac{1}{\epsilon}}(|\dev(T)|^{p+1})\,v\,\di x\\[1ex]
&+&\int_{\partial\Omega}g_{\theta}\,v\,\di x
\end{eqnarray*}
for all $v\in H^1(\Omega;\R)$ and almost all $t\in(0,T)$.
\end{de}
\begin{tw}$\mathrm{(Main\, Theorem)}$\\[1ex]
\label{tw:1.2}
Suppose that the given data have the regularity required in (\ref{eq:1.4})-(\ref{eq:1.6}). Moreover, let $$|\dev(\D(\ve(u_0)-\ve^{p,0})|^{p-1}\dev(\D(\ve(u_0)-\ve^{p,0})\in L^2(\Omega;\SS).$$ Then, for all $T>0$ the system (\ref{eq:1.1}) with the boundary condition (\ref{eq:1.2}) and the initial condition (\ref{eq:1.3}) possesses a solution in the sense of Definition \ref{de:1.1}.
\end{tw}
The initial assumption in Theorem \ref{tw:1.2} means that for $t=0$ the argument of the
maximal monotone operator belongs to the domain of $|\dev(T)|^{p-1}\dev(T)$ (for more information we refer
to \cite{26}, Definition 2.4 and Theorem 2.5).\\
Our Theorem \ref{tw:1.2} will be proved in the next three sections. First, we use the Yosida approximation to the maximal monotone inelastic constitutive equation. Then we pass to the limit to obtain a solution in the sense of Definition \ref{de:1.1}.
\section{Transformation to a homogeneous boundary-value problem with respect to the temperature}


\renewcommand{\theequation}{\thesection.\arabic{equation}}
\setcounter{equation}{0}
Let us consider the following linear parabolic system
\begin{eqnarray}
\label{eq:2.4}
\tilde{\theta}_t(x,t)-\Delta\tilde{\theta}(x,t)=0
\end{eqnarray}
with boundary-initial conditions
\begin{eqnarray}
\label{eq:2.5}
\frac{\partial\,\tilde{\theta}}{\partial\,n}(x,t)&=&g_{\theta}(x,t)\quad \textrm{ for}\quad x\in \partial \Omega \quad\textrm{and}\quad  t\geq 0\,,\nn\\[1ex]
\tilde{\theta}(x,0)&=&\theta_0(x)\quad\,\,\,\, \textrm{ for}\quad x\in \Omega.
\end{eqnarray}
Assuming that  $g_{\theta}\in H^{1}(0,T;L^2(\partial\Omega;\R))$ and $\theta_0\in H^1(\Omega;\R)$ we conclude that the system (\ref{eq:2.4}) possesses a solution $\tilde\theta\in L^{\infty}(0,T;H^1(\Omega;\R))$ such that $\tilde{\theta}_t\in L^2(0,T;L^2(\Omega;\R))$. Additionally the following estimate 
\begin{eqnarray}
\label{eq:2.6}
\|\tilde{\theta}_t\|_{L^{2}(0,T;L^2(\Omega;\R))}+\|\tilde{\theta}\|_{L^{\infty}(0,T;H^1(\Omega;\R))}\leq D \,\Big(\|g_{\theta}\|_{H^{1}(0,T;L^2(\partial\Omega;\R))} + \|\theta_0\|_{H^1(\Omega;\R)}\Big)
\end{eqnarray}
holds. Finally, we define $\theta=\hat{\theta}-\tilde{\theta}$. Notice that to find a solution $(u,\hat{\theta})$ to the problem (\ref{eq:1.1})-(\ref{eq:1.3}) we need to find a solution $(u,\theta)$ of the following problem
\begin{eqnarray}
\label{eq:2.7}
\mathrm{div}_x \,\sigma&=&-F-\mathrm{div}\,\D (\ve(u_t))\,,\nn\\
\sigma&=&\D(\ve(u)-\ve^p)-f\big(T_{\frac{1}{\epsilon}}(\theta+\tilde{\theta})\big)\id\,,\nn\\[1ex]
\ve^{p}_{t}&=& |\dev(T)|^{p-1}\,\dev(T)\,,\\[1ex]
T&=&\D(\ve(u)-\ve^p)\,,\nn\\[1ex]
\theta_t-\Delta\theta+f\big( T_{\frac{1}{\epsilon}}(\theta+\tilde{\theta})\big)\mathrm{div}u_t&=& T_{\frac{1}{\epsilon}}(|\dev(T)|^{p+1})\nn
\end{eqnarray}
with the initial-boundary conditions 
\begin{eqnarray}
\label{eq:2.8}
u_{|_{\partial\Omega}}=g_D\,,&&
\frac{\partial\,\theta}{\partial\,n}_{|_{\partial\Omega}}=0\,,\nn\\[1ex]
\theta(0)=\hat{\theta}_0-\tilde{\theta}_0=\theta_{0}\,,&&
u(0)=u_0\,,\quad
\ve^p(0)=\ve^{p,0}\,.
\end{eqnarray}
\section{Yosida approximation}
Notice that the inelastic constitutive equation $(\ref{eq:2.7})_3$ is a maximal monotone vector field such that $G(0)=0$, where $G(T)=|\dev(T)|^{p-1}\dev(T)$. Moreover, there exists a positive, differentiable convex function $M : \SS\rightarrow\SS$ such that $\nabla_T M=G$. The main idea to prove an existence result for system (\ref{eq:2.7}) is based on the, so called, partial Yosida approximation. This means that we will use the Yosida approximation of the maximal monotone vector field $\nabla M$. For $\lambda>0$ let us define the function
\begin{eqnarray*}
M_{\lambda}(z)=\inf_{w\in\R^6}\frac{1}{2\lambda}\{|z-w|^2+M(w)\}\textrm{.}
\end{eqnarray*}
Hence, $M_{\lambda}$ is subquadratic, nonnegative, the gradient $\nabla M_{\lambda}$ is a global Lipschitz map and is the Yosida approximation of $\nabla M$ - for more information we refer to \cite{23} and \cite{24}. Using the function $M_{\lambda}$ we define a sequence of approximate problems
\renewcommand{\theequation}{\thesection.\arabic{equation}}
\setcounter{equation}{0}%
\begin{eqnarray}
\label{eq:3.1}
\mathrm{div}_x \,\sigma^{\lambda}&=&-F-\mathrm{div}\,\D(\ve^{\lambda}(u_t))\,,\nn\\
\sigma^{\lambda}&=&\D(\ve(u^{\lambda})-\ve^{p,\lambda})- f\big(T_{\frac{1}{\epsilon}}(\theta^{\lambda}+\tilde{\theta})\big)\id\,,\nn\\[1ex]
\ve^{p,\lambda}_{t}&=&\nabla M_{\lambda}(\dev(T^{\lambda}))\,,\\[1ex]
T^{\lambda}&=&\D(\ve(u^{\lambda})-\ve^{p,\lambda})\,,\nn\\[1ex]
\theta_t^{\lambda}-\Delta\theta^{\lambda}+ f\big(T_{\frac{1}{\epsilon}}(\theta^{\lambda}+\tilde{\theta})\big)\mathrm{div}\,u_t^{\lambda}&=& T_{\frac{1}{\epsilon}}(\dev(T^{\lambda})\,\ve^{p,\lambda}_t)\,.\nn
\end{eqnarray}
The problem (\ref{eq:3.1}) is considered with the boundary conditions
\begin{eqnarray}
\label{eq:3.2}
u^{\lambda}(x,t)&=&g_{D}(x,t)\qquad x\in \partial\Omega \textrm{,}\quad t\geq 0\textrm{,}\nn\\[1ex]
\frac{\partial\,\theta^{\lambda}}{\partial\,n}(x,t)&=&0\qquad\qquad\,\,\, x\in \partial \Omega, \quad  t\geq 0
\end{eqnarray}
and with the initial conditions
\begin{eqnarray}
\label{eq:3.3}
u^{\lambda}(x,0)=u_0(x),\quad\ve^{p,\lambda}(x,0)=\ve^{p,0}(x),\quad \theta^{\lambda}(x,0)=\theta_0(x).
\end{eqnarray}
To prove existence result for all $\lambda>0$ we are going to define the operator $\mathcal{R}$ acting from the Banach space $L^{2}(0,T;H^1(\Omega;\R))$ into the same space. Next, using to this operator the Banach Fixed Point Theorem we find a solution of system (\ref{eq:3.1}) with initial-boundary conditions (\ref{eq:3.2}) and (\ref{eq:3.3}).\\[1ex]
Fix $\lambda>0$. Let us set $\theta^{\star}\in L^{2}(0,T;H^1(\Omega;\R))$ and consider the visco-elasticity problem:
\begin{eqnarray}
\label{eq:3.4}
\mathrm{div}_x \sigma&=&-F-\mathrm{div}\,\D (\ve(u_t))\,,\nn\\
\sigma&=&\D(\ve(u)-\ve^{p})- f\big(T_{\frac{1}{\epsilon}}(\theta^{\star}+\tilde{\theta})\big)\id\,,\nn\\[1ex]
\ve^{p}_{t}&=&\nabla M_{\lambda}(\dev(T))\,,\\[1ex]
T&=&\D(\ve(u)-\ve^{p})\,,\nn\\[1ex]
u_{|_{\partial\Omega}}&=&g_D\,,\nn\\[1ex]
u(0)&=&u_0\,,\nn\\[1ex]
\ve^p(0)&=&\ve^{p,0}\nn\,.
\end{eqnarray}
Notice that in the system (\ref{eq:3.4}) we dropped the subscript $\lambda>0$. In the next part of this section, we will drop the subscript $\lambda$ and write $u$, $\theta$ and $\ve^{p}$
instead of $u^{\lambda}$, $\theta^{\lambda}$ and $\ve^{p,\lambda}$.
\begin{lem}
\label{lem:3.1}
Let us assume that the given data satisfy all requirements of Theorem \ref{tw:1.2} and $\theta^{\star}\in L^{2}(0,T;H^1(\Omega;\R))$. Then, there exists a global in time unique solution $(u,\ve^p)$ of the system (\ref{eq:3.4}) such that 
$$(u,\ve^p)\in H^1(0,T;H^1(\Omega;\R^3))\times W^{1,\infty}(0,T;L^2(\Omega;\SS)).$$
\end{lem}
{\bf\em Proof:}\hspace{2ex} The idea of the proof is a fixed-point argument. Let $\ve^{\ast}\in L^{\infty}(0,T;L^2(\Omega;\S))$ and let us consider the equation
\begin{equation}
\label{eq:3.5}
\ve^p(x,t)=\ve^{p,0}(x)+\int_0^t \nabla M_{\lambda}\Big(\dev\big(\D(\ve^{\ast}(\tau)-\ve^p(\tau))\big)\Big)\di \tau\,.
\end{equation}
From the theory of differential equations in Banach spaces ($\nabla M_{\lambda}$ is global Lipschitz) it follows that the equation (\ref{eq:3.5}) possesses a global in time, unique solution\\ $\ve^p\in W^{1,\infty}(0,T;L^2(\Omega;\SS))$. Next, consider the system of equations
\begin{eqnarray}
\label{eq:3.6}
\mathrm{div}_x \sigma&=&-F-\mathrm{div}\,\D( \ve(u_t))\,,\nn\\[1ex]
\sigma&=&\D(\ve(u)-\ve^{p})- f\big(T_{\frac{1}{\epsilon}}(\theta^{\star}+\tilde{\theta})\big)\id\,,\\[1ex]
u(0)=u_0\,,&&
u_{|_{\partial\Omega}}=g_D\nn\,,
\end{eqnarray}
where $\ve^p$ satisfies (\ref{eq:3.5}). Notice that the term $f\big(T_{\frac{1}{\epsilon}}(\theta^{\star}+\tilde{\theta}))$ is bounded, because the argument of the continuous function $f$ is bounded. Therefore the system (\ref{eq:3.6}) has a unique solution $u\in H^1(0,T;H^1(\Omega;\R^3))$. This way, we defined an operator
\begin{equation}
\label{eq:3.7}
P:L^{\infty}(0,T;L^2(\Omega;\S))\rightarrow L^{\infty}(0,T;L^2(\Omega;\S))
\end{equation}
such that $P(\ve^{\ast})=\frac{1}{2}(\nabla\,u+\nabla^Tu)$. We will show that $P$ is a contraction.\\[1ex]
Let us denote by $\ve^p_1(t)$ and $\ve^p_2(t)$ solutions of (\ref{eq:3.5}) with the input functions\\ $\ve^{\ast}_1,\,\ve^{\ast}_2\in L^{\infty}(0,T;L^2(\Omega;\S))$, respectively. The difference $P(\ve^{\ast}_1)-P(\ve^{\ast}_2)$ satisfies
\begin{eqnarray}
\label{eq:3.8}
-\mathrm{div}\,\big(\D(P(\ve^{\ast}_1)-P(\ve^{\ast}_2)+P(\ve^{\ast}_{1})_t-P(\ve^{\ast}_{2})_t)\big) &=& -\mathrm{div}\,\big(\D(\ve^p_1-\ve^p_2)\big)
\end{eqnarray} 
and the following inequality
\begin{eqnarray}
\label{eq:3.9}
\|P(\ve^{\ast}_1(t))-P(\ve^{\ast}_2(t))\|_{L^2(\Omega;\S)}\leq C\|\ve^p_1(t)-\ve^p_2(t)\|_{L^2(\Omega;\SS)}
\end{eqnarray}
holds, where the positive constant $C$ does not depend on these input functions, $\theta^{\star}$ and is independent of $t$. Using equation (\ref{eq:3.5}) it is not difficult to obtain the following inequality 
\begin{eqnarray}
\label{eq:3.10}
\|\ve^p_1(t)-\ve^p_2(t)\|_{L^2(\Omega;\SS)}\leq \tilde{C}\|\ve^{\ast}_1(t)-\ve^{\ast}_2(t)\|_{L^2(\Omega;\S)},
\end{eqnarray}
where $\tilde{C}$ does not depend on $t$ (it depends only on the Lipschitz constant and on time $T$). Having these two inequalities we easily see that the operator $P$ is a contraction- for more details we refer to \cite{26}..\\$\mbox{}$ \hfill $\Box$
\begin{tw}
Suppose that the given data satisfy all requirements of Theorem \ref{tw:1.2}. 
For all $\lambda>0$ the system (\ref{eq:3.1}) with initial-boundary conditions (\ref{eq:3.2}) and (\ref{eq:3.3}) possesses unique, global in time solution $(u,\ve^p,\theta)$ such that 
$$(u,\ve^p)\in H^1(0,T;H^1(\Omega;\R^3))\times W^{1,\infty}(0,T;L^2(\Omega;\SS))\,,$$
$$\theta\in L^{\infty}(0,T;H^1(\Omega;\R))\,\,\mathrm{and}\,\, \theta_t\in L^{2}(0,T;L^2(\Omega;\R)).$$
\end{tw}
{\bf\em Proof:}\hspace{2ex} Fix $\theta^{\star}\in L^{2}(0,T;H^1(\Omega;\R))$. Lemma \ref{lem:3.1} implies that there exists solution $(u,\ve^p)\in H^1(0,T;H^1(\Omega;\R^3))\times W^{1,\infty}(0,T;L^2(\Omega;\SS))$ of the system (\ref{eq:3.4}). Inserting this solution to heat equation $(\ref{eq:3.1})_5$ we obtain a solution $\theta\in L^{\infty}(0,T;H^1(\Omega;\R))$ such that $\theta_t\in L^{2}(0,T;L^2(\Omega;\R))$. Hence, we have defined an operator
$$L^{2}(0,T;H^1(\Omega;\R))\owns\theta^{\star}\rightarrow\mathcal{R}(\theta^{\star})=\theta\in L^{2}(0,T;H^1(\Omega;\R)).$$
Next we prove that $\mathcal{R}$ is a contraction. Let $(u_1,\ve^{p}_1)$ and $(u_2,\ve^{p}_2)$ be the solutions to (\ref{eq:3.4}) for $\theta^{\star}_1$ and $\theta^{\star}_2\in L^{2}(0,T;H^1(\Omega;\R))$, respectively. Let us mark by $\bar{u}=u_1-u_2$ and $\ve^p=\ve^p_1-\ve^p_2$. Then $\ve^p(0)=0$ and $\bar{u}(0)=0$. The weak formulation of the system (\ref{eq:3.4}) yields
\begin{eqnarray}
\label{eq:3.11}
\int_{\Omega}\D(\ve(\bar{u})+\ve(\bar{u}_t)-\ve^p)\ve(v)\,\di x = \int_{\Omega} \Big(f\big(T_{\frac{1}{\epsilon}}(\theta^{\star}_1+\tilde{\theta})\big)- f(T_{\frac{1}{\epsilon}}(\theta^{\star}_2+\tilde{\theta})\big)\Big) \mathrm{div}\,v\,\di x 
\end{eqnarray} 
for all $v\in H_0^1(\Omega,\R^3)$. Putting $v=\bar{u}_t$ in (\ref{eq:3.11}) we obtain
\begin{eqnarray}
\label{eq:3.12}
&&\frac{1}{2}\frac{d}{dt}\Big(\int_{\Omega}\D(\ve(\bar{u})-\ve^p)(\ve(\bar{u})-\ve^p)\,\di x\Big) + \int_{\Omega}\D\ve(\bar{u}_t)\ve(\bar{u}_t)\,\di x=
-\int\limits_{\Omega}\D(\ve(\bar{u})-\ve^p)\ve^p_t\,\di x\nn\\[1ex]
&+&\int_{\Omega} \Big(f\big(T_{\frac{1}{\epsilon}}(\theta^{\star}_1+\tilde{\theta})\big)- f(T_{\frac{1}{\epsilon}}(\theta^{\star}_2+\tilde{\theta})\big)\Big)\mathrm{div}\,\bar{u}_t\,\di x 
\end{eqnarray}
Notice that the difference $f\big(T_{\frac{1}{\epsilon}}(\theta^{\star}_1+\tilde{\theta})\big)- f(T_{\frac{1}{\epsilon}}(\theta^{\star}_2+\tilde{\theta})\big)$ is bounded. Using monotonicity of $\nabla M_{\lambda}$ and integrating with respect to time we get  
\begin{eqnarray}
\label{eq:3.13}
\int_0^T\int_{\Omega}\D\ve(\bar{u}_t)\ve(\bar{u}_t)\,\di x\,\di t \leq T\,C
\end{eqnarray}
and the constant $C>0$ does not depend on $\theta^{\star}_1$ and $\theta^{\star}_2$ (it depends on $\epsilon$ and $\Omega$ only). Additionally, from the weak formulation of $(\ref{eq:3.1})_5$ we have
\begin{eqnarray}
\label{eq:3.14}
&&\frac{1}{2}\frac{d}{dt}\Big(\int_{\Omega}|\theta_1-\theta_2|^2\,\di x\Big) + \int_{\Omega}|\nabla(\theta_1-\theta_2)|^2\,\di x=\nn\\[1ex]
&&\int_{\Omega} \Big(f\big(T_{\frac{1}{\epsilon}}(\theta_1+\tilde{\theta})\big)\mathrm{div}\,u^1_t- f\big(T_{\frac{1}{\epsilon}}(\theta_2+\tilde{\theta}))\mathrm{div}\,u_t^2\Big)(\theta_1-\theta_2)\di x\\[1ex]
&+&\int_{\Omega}\Big(T_{\frac{1}{\epsilon}}(\dev(T_1)\,\ve^{p,1}_t) -T_{\frac{1}{\epsilon}}(\dev(T_2)\,\ve^{p,2}_t)\Big) (\theta_1-\theta_2)\di x\nn\,,
\end{eqnarray}
where $\theta^1$ and $\theta^2$ are solutions of $(\ref{eq:3.1})_5$ with the input functions  $(u_1,\ve^{p}_1)$ and $(u_2,\ve^{p}_2)$, respectively.
Using the Cauchy-Schwarz inequality with small weight and integrating with respect to time we get
\begin{eqnarray}
\label{eq:3.15}
&&\int_{\Omega}|\theta_1-\theta_2|^2\,\di x + \int^T_0\int_{\Omega}|\nabla(\theta_1-\theta_2)|^2\,\di x\,\di t\leq
C(\nu)\int^T_0\int_{\Omega} |\ve(u^1_t-u^2_t)|^2\di x\,\di t\nn\\[1ex] &&\nu\int^T_0\int_{\Omega}|\theta_1-\theta_2|^2\di x\,\di t +T\,\tilde{C}\,,
\end{eqnarray}
where the constant $\tilde{C}>0$ does not depend on $\theta^{\star}_1$ and $\theta^{\star}_2$. Inserting (\ref{eq:3.13}) into (\ref{eq:3.15}) and choosing $\nu>0$ sufficient small we receive
\begin{eqnarray}
\label{eq:3.16}
\int_{\Omega}|\theta_1-\theta_2|^2\,\di x + \int^T_0\int_{\Omega}|\nabla(\theta_1-\theta_2)|^2\,\di x\,\di t\leq T\,D\int^T_0\|\theta^{\star}_1-\theta^{\star}_2\|_{H^1(\Omega;\R)}^2\,\di t,
\end{eqnarray}
where the constant $D>0$ does not depend on $\theta^{\star}_1$ and $\theta^{\star}_2$ (it depends on $\epsilon$ and $\Omega$ only). Moreover, notice that $D$ is also independent on the initial data $\ve^{p,0}$, $u_0$ and $\theta_0$, thus the inequality (\ref{eq:3.16}) implies that the operator $\mathcal{R}$ is a contraction. $\mbox{}$ \hfill $\Box$


\section{Passing to the limit in Yosida approximation}
This section is devoted to prove some a priori estimates for the sequence of approximate solutions $\{(u^{\lambda},\theta^{\lambda},\ve^{p,\lambda})\}_{\lambda>0}$ and pass to the limit $\lambda\rightarrow 0^+$. In this section we returned to the subscription $\lambda$.\\[1ex]
{\bf Remark:}\hspace{2ex} The construction of the operator $\mathcal{R}$ yields that $\theta_t^{\lambda}\in L^2(0,T;L^2(\Omega;\R))$ and we can conclude that the following elliptic problem
\begin{eqnarray*}
\mathrm{div}\Big(\D\big(\ve(u^{\lambda}_t)-\ve^{p,\lambda}_t\big)- f'\big(T_{\frac{1}{\epsilon}}(\theta^{\lambda}+\tilde{\theta})\big) \nabla T_{\frac{1}{\epsilon}}(\theta^{\lambda}+\tilde{\theta}) (\theta^{\lambda}_t+\tilde{\theta}_t)\id\Big)&=& -F_t-\mathrm{div}\,\D(\ve(u_{tt}^{\lambda}))\,,\\[1ex]
u^{\lambda}_t\mathrm{}_{|_{\Omega}}&=&g_{D,t}
\end{eqnarray*}
possesses the solution $u_t^{\lambda}\in H^1(0,T;H^1(\Omega;\R^3))$ ($ T_{\frac{1}{\epsilon}}$ is a Lipschitz function and the chain rule for the derivation of $ T_{\frac{1}{\epsilon}}$ holds true).\\[2ex]
Let us start with prove some estimates for the approximation sequence.
\renewcommand{\theequation}{\thesection.\arabic{equation}}
\setcounter{equation}{0}%
\begin{tw} $\mathrm{(energy\, estimate)}$\\
\label{tw:4.1}
Assume that the given data satisfy all requirements of Theorem \ref{tw:1.2}.
Then there exists a positive constant $C(T)$ (not depending on $\lambda$) such that the following inequality holds
\begin{eqnarray*}
\|T^{\lambda}(t)\|^2_{L^2(\Omega;\S)} + \int_0^t\|\ve^{\lambda}_t(\tau)\|^2_{L^2(\Omega;\S)}\,\di\tau\leq C(T)
\end{eqnarray*}
for $t\leq T$.
\end{tw}
{\bf Proof:}\hspace{2ex} Calculate the time derivative 
\begin{eqnarray}
\label{eq:4.1}
&&\frac{d}{dt}\Big(\frac{1}{2}\int_{\Omega}\D(\ve^{\lambda}-\ve^{p,\lambda}) (\ve^{\lambda}-\ve^{p,\lambda})\,\di x\Big)= \int_{\Omega}\D(\ve^{\lambda}-\ve^{p,\lambda})
(\ve^{\lambda}_t-\ve_t^{p,\lambda})\,\di x\nn\\[1ex]
&=&\int_{\Omega}\sigma^{\lambda}\,\ve^{\lambda}_t\,\di x 
+\int_{\Omega} f\big(T_{\frac{1}{\epsilon}}(\theta^{\lambda}+\tilde{\theta}))\,\mathrm{div}\,u^{\lambda}_t\,\di x
-\int_{\Omega}\dev(T^{\lambda})\,\ve^{p,\lambda}_t\,\di x\,.\qquad\
\end{eqnarray}
The last integral on the right-hand side of (\ref{eq:4.1}) is non-negative. Integrating with respect to time we obtain
\begin{eqnarray}
\label{eq:4.2}
&&\int_{\Omega}\D(\ve^{\lambda}-\ve^{p,\lambda}) (\ve^{\lambda}-\ve^{p,\lambda})\,\di x
\leq \int_{\Omega}\D(\ve^{\lambda}(0)-\ve^{p,0}) (\ve^{\lambda}(0)-\ve^{p,0})\,\di x\nn\\[1ex] &+&\int_0^t\int_{\Omega}\sigma^{\lambda}\,\ve^{\lambda}_t\,\di x\,\di\tau
+\int_0^t\int_{\Omega} f\big(T_{\frac{1}{\epsilon}}(\theta^{\lambda}+\tilde{\theta})\big)\,\mathrm{div}\,u^{\lambda}_t\,\di x\,\di\tau\,.
\end{eqnarray}
Integrating by parts in the second term on the right hand-side of (\ref{eq:4.2}) and  using the equation $(\ref{eq:3.1})_{1}$ we get
\begin{eqnarray}
\label{eq:4.3}
&&\int_{\Omega}\D(\ve^{\lambda}-\ve^{p,\lambda}) (\ve^{\lambda}-\ve^{p,\lambda})\,\di x + \int_0^t\int_{\Omega}\D(\ve^{\lambda}_{t})\ve^{\lambda}_t\,\di x\, \di\tau\nn\\[1ex] 
&\leq&\int_{\Omega}\D(\ve^{\lambda}(0)-\ve^{p,0}) (\ve^{\lambda}(0)-\ve^{p,0})\,\di x 
+\int_0^t\int_{\Omega}F\,u^{\lambda}_t\,\di x\,\di\tau\nn\\[1ex]
&+& \int_0^t\int_{\partial\Omega}\big(\sigma^{\lambda}+\D\big(\ve(u^{\lambda}_t)\big)\big)n\,g_{D,t}\,\di S\,\di\tau
+\int_0^t\int_{\Omega} f\big(T_{\frac{1}{\epsilon}}(\theta^{\lambda}+\tilde{\theta})\big)\,\mathrm{div}\,u^{\lambda}_t\,\di x\,\di\tau\,.
\end{eqnarray}
The boundary integral appearing on the right-hand side of inequality (\ref{eq:4.3}) can be estimated using the continuity of the trace operator in the space $L^2(\mathrm{div})$ [27]
\begin{eqnarray}
\label{eq:4.4}
\lefteqn{\int_0^t\int_{\partial\Omega}\big(\sigma^{\lambda}+\D\big(\ve(u^{\lambda}_t)\big)\big)n\,g_{D,t}\,\di S\,\di\tau}\nn\\[1ex] &\leq& c\big(\int_0^t\|\big(\sigma^{\lambda}+\D\big(\ve(u^{\lambda}_t)\big)\big)n\|_{H^{-\frac{1}{2}}(\partial\Omega;\R^3)} \|g_{D,t}\|_{H^{\frac{1}{2}}(\partial\Omega;\R^3)}\,\di\tau\big)\nn\\[1ex] &\leq& c\Big(\int_0^t\big(\|\sigma^{\lambda}+\D\big(\ve(u^{\lambda}_t)\big)\|_{L^2(\Omega;\S)}+ \|\textrm{div}\,\big(\sigma^{\lambda}+\D\big(\ve(u^{\lambda}_t\big)\big)\|_{L^2(\Omega;\R^3)}\big) \|g_{D,t}\|_{H^{\frac{1}{2}}(\partial\Omega;\R^3)}\,\di\tau\Big)\nn\\[1ex]
&\leq& c\Big(\int_0^t\big(\|\sigma^{\lambda}+\D\big(\ve(u^{\lambda}_t)\big)\|_{L^2(\Omega;\S)}+ \|F\|_{L^2(\Omega;\R^3)}\big) \|g_{D,t}\|_{H^{\frac{1}{2}}(\partial\Omega;\R^3)}\,\di\tau\Big)\nn\\[1ex]
&\leq& c\Big(\int_0^t\big(\|T^{\lambda}\|_{L^2(\Omega;\S)}+ \|F\|_{L^2(\Omega;\R^3)}\big) \|g_{D,t}\|_{H^{\frac{1}{2}}(\partial\Omega;\R^3)}\,\di\tau\nn\\[1ex]
&+&\int_0^t\|\ve^{\lambda}_t\|_{L^2(\Omega;\S)} \|g_{D,t}\|_{H^{\frac{1}{2}}(\partial\Omega;\R^3)}\,\di\tau+1\Big).\qquad
\end{eqnarray}
Putting (\ref{eq:4.4}) into (\ref{eq:4.3}) we receive
\begin{eqnarray}
\label{eq:4.5}
&&\int_{\Omega}\D(\ve^{\lambda}-\ve^{p,\lambda}) (\ve^{\lambda}-\ve^{p,\lambda})\,\di x + \int_0^t\int_{\Omega}\D(\ve^{\lambda}_{t})\ve^{\lambda}_t\,\di x\, \di\tau\,\,\leq\,\,\tilde{C}(T)\nn\\[1ex] 
&\leq& \int_{\Omega}\D(\ve(u_0)-\ve^{p,0}) (\ve(u_0)-\ve^{p,0})\,\di x +
\nu\int_0^t\int_{\Omega}|\ve^{\lambda}_{t}|^2\,\di x\,\di\tau+\nu \int_0^t\int_{\Omega}|T^{\lambda}|^2\,\di x\,\di\tau\,,\qquad
\end{eqnarray}
where $\nu>0$ is any positive constant and $\tilde{C}(T)$ does not depend on $\lambda$. Finally, choosing in (\ref{eq:4.5}) $\nu>0$ sufficient small we finish the proof.$\mbox{}$ \hfill $\Box$\\[1ex]
{\bf Remark:}\hspace{2ex} The Theorem \ref{tw:4.1} implies that the sequence $\{u^{\lambda}_t\}_{\lambda>0}$ is bounded in the space $L^2(0,T;H^1(\Omega;\R^3))$. Let us consider the equation $(\ref{eq:3.1})_5$ 
\begin{eqnarray}
\label{eq:4.8}
\theta_t^{\lambda}-\Delta\theta^{\lambda}&=& T_{\frac{1}{\epsilon}}(\dev(T^{\lambda})\,\ve^{p,\lambda}_t) - f\big(T_{\frac{1}{\epsilon}}(\theta^{\lambda}+\tilde{\theta})\big)\mathrm{div}\,u_t^{\lambda}\,.
\end{eqnarray}
Notice that the right hand-side of (\ref{eq:4.8}) is bounded independently on $\lambda$ in the space $L^2(0,T;L^2(\Omega;\R))$. From the standard theory for heat equation we immediately conclude that the sequences $\{\theta^{\lambda}\}_{\lambda>0}$ and $\{\theta^{\lambda}_t\}_{\lambda>0}$ are bounded in the space $L^{\infty}(0,T;H^1(\Omega;\R))$ and $L^2(0,T;L^2(\Omega;\R))$, respectively.\\[2ex]
The next step in our existence theory is an estimate for time derivatives of the approximate sequence. The proof of the theorem below works for gradient flows only. The main idea is similar to that from Theorem 2.3 in \cite{28}. 
\begin{tw}\label{tw:4.2}$\mathrm{(L^2(L^2)\,estimates\, for\, time\, derivatives)}$\\[0.5ex]
Assume that the given data satisfy all requirements of Theorem \ref{tw:1.2}. Then, for all $T>0$ and $t\leq T$ the solution of approximate problem (\ref{eq:3.1}) - (\ref{eq:3.3}) satisfies the following inequality
\begin{eqnarray}
\label{eq:4.9}
\int\limits_0^t\|T^{\lambda}_t(\tau)\|^2_{L^2(\Omega;\S)}\,\di\tau+
\|\ve^{\lambda}_t(\tau)\|^2_{L^2(\Omega;\S)}\leq C(T)\,\textrm{,}
\end{eqnarray}
where the constant $C(T)$ does not depend on $\lambda$.
\end{tw}
{\em Proof:}\hspace{2ex} Let us fix $T>0$. Compute the derivative 
\begin{eqnarray}
\label{eq:4.10}
\lefteqn{\frac{d}{dt}\Big(\int_{\Omega}M_{\lambda}(\dev(T^{\lambda}))\,\di x\Big)=
\int_{\Omega}\nabla M_{\lambda}(\dev(T^{\lambda}))(\dev(T^{\lambda}_t))\,\di x}\nn\\[1ex]
&=&\int_{\Omega}\dev\big(\D(\ve(u_t^{\lambda})-\ve^{p,\lambda}_t)\big)\ve^{p,\lambda}_t\,\di x=\int_{\Omega}T^{\lambda}_t\,\ve^{p,\lambda}_t\,\di x
=-\int_{\Omega}\D^{-1}T^{\lambda}_tT^{\lambda}_t\,\di x\nn\\[1ex]
&+&\int_{\Omega}T^{\lambda}_t\ve^{\lambda}_t\,\di x
= -\int_{\Omega}\D^{-1}T^{\lambda}_tT^{\lambda}_t\,\di x+\int_{\Omega}\sigma^{\lambda}_t\,\ve^{\lambda}_t\,\di x\nn\\[1ex]
&+&\int_{\Omega} f'\big(T_{\frac{1}{\epsilon}}(\theta^{\lambda}+\tilde{\theta})\big) \nabla T_{\frac{1}{\epsilon}}(\theta^{\lambda}+\tilde{\theta})\,(\theta^{\lambda}_t+ \tilde{\theta}_t)\,\mathrm{div}\,u^{\lambda}_t\,\di x\,.
\end{eqnarray}
Notice that in the formula (\ref{eq:4.10}) we use the following informations: $T_{\frac{1}{\epsilon}}$ is a Lipschitz function, $T_{\frac{1}{\epsilon}}\big(\theta^{\lambda}+\tilde{\theta}\big)$ belongs to $H^1(\Omega)$ and the chain rule for derivation of $T_{\frac{1}{\epsilon}}\big(\theta^{\lambda}+\tilde{\theta}\big)$ holds true (see e.g. \cite{29,30}).
Integrating the equality (\ref{eq:4.10}) with respect to time, we get
\begin{eqnarray}
\label{eq:4.11}
\lefteqn{\int_{\Omega}M_{\lambda}(\dev(T^{\lambda}(t)))\,\di x + \int_0^t\int_{\Omega}\D^{-1}T^{\lambda}_t T^{\lambda}_t\,\di x\,\di\tau+\frac{1}{2}\int_{\Omega}\D(\ve_{t}^{\lambda})\ve_{t}^{\lambda} \,\di x}\nn\\[1ex] 
&=&\int_{\Omega}M_{\lambda}(\dev(T^{\lambda}(0)))\,\di x+ \frac{1}{2}\int_{\Omega}\D(\ve_{t}^{\lambda}(0))\ve_{t}^{\lambda}(0) \,\di x+\int_0^t\int_{\Omega}\big(\sigma^{\lambda}_t+\D(\ve_{tt}^{\lambda})\big)\,\ve^{\lambda}_t\,\di x\,\di\tau\nn\\[1ex]
&+&\int_0^t\int_{\Omega}f'\big(T_{\frac{1}{\epsilon}}(\theta^{\lambda}+\tilde{\theta})\big)\nabla T_{\frac{1}{\epsilon}}\big(\theta^{\lambda}+\tilde{\theta}\big)\,(\theta^{\lambda}_t+ \tilde{\theta}_t)\,\mathrm{div}\,u^{\lambda}_t\,\di x\,.
\end{eqnarray}
The function $u_t(0)$ is the solution of the following linear elliptic problem
\begin{eqnarray}
\label{eq:4.12}
\mathrm{div}\,\D\big(\ve(u_t^{\lambda}(0))\big)&=& -F(0)-\mathrm{div}\Big(\D(\ve(u_0)-\ve^{p,0})- f\big(T_{\frac{1}{\epsilon}}(\theta_0)\big)\id\Big)\,,\nn\\[1ex]
u^{\lambda}(x,0)=u_0(x), &&\quad u_t(0)_{|_{\partial\Omega}}=g_{D,t}(0)\,.
\end{eqnarray}
Therefore, we obtain inequality
\begin{eqnarray}
\label{eq:4.13}
\|u^{\lambda}_t(0)\|^2_{H^1(\Omega;\R^3)}&\leq& \underline{C}\Big(\|u_0\|_{H^1(\Omega;\R^3)}^2+\|F(0)\|_{L^2(\Omega;\R^3)}^2 +\|\ve^{p,0}\|_{L^2(\Omega;\SS)}^2\nn\\[1ex]
&+&\|f\big(T_{\frac{1}{\epsilon}}(\theta_0)\big)\|_{L^2(\Omega;\R)}^2
+\|g_{D,t}(0)\|_{H^{\frac{1}{2}}(\partial\Omega;\R^3)}^2\Big)\,.
\end{eqnarray}
The inequality (\ref{eq:4.13}) implies that $\|u^{\lambda}_t(0)\|_{H^1(\Omega;\R^3)}$ is bounded (independently of $\lambda$). The function $T_{\frac{1}{\epsilon}}$ is Lipschitz, hence $\nabla T_{\frac{1}{\epsilon}}$ is bounded. Moreover, the remark before the Theorem \ref{tw:4.2} implies that the last term on the right hand-side of (\ref{eq:4.11}) is bounded. Now we integrate by parts in the third term on the right hand-side of (\ref{eq:4.11}) to get  
\begin{eqnarray}
\label{eq:4.14}
&&\int_0^t\int_{\Omega}\big(\sigma^{\lambda}_t+\D(\ve_{tt}^{\lambda})\big)\,\ve^{\lambda}_t\,\di x\,\di\tau=
\int_0^t\int_{\Omega}F_t\,u^{\lambda}_t\,\di x\,\di\tau + \int_0^t\int_{\partial\Omega}\big(\sigma^{\lambda}_t+\D(\ve_{tt}^{\lambda})\big)n\,g_{D,t}\,\di S\,\di\tau\nn\\[1ex]
&=&\int_0^t\int_{\Omega}F_t\,u^{\lambda}_t\,\di x\,\di\tau - \int_0^t\int_{\partial\Omega}\big(\sigma^{\lambda}+\D(\ve_{t}^{\lambda})\big)n\,g_{D,tt}\,\di S\,\di\tau\nn\\[1ex] &+&\int_{\partial\Omega}\big(\sigma^{\lambda}(t)+\D(\ve_{t}^{\lambda}(t))\big)n\,g_{D,t}(t)\,\di S
-\int_{\partial\Omega}\big(\sigma^{\lambda}(0)+\D(\ve_{t}^{\lambda}(0))\big)n\,g_{D,t}(0)\,\di S\,.
\end{eqnarray}
The boundary integrals are estimated using the continuity of the trace operator in the space $L^2(\mathrm{div})$ (in the same way as in (\ref{eq:4.4})) hence we get
\begin{eqnarray}
\label{eq:4.15}
\int_0^t\int_{\Omega}\big(\sigma^{\lambda}_t+\D(\ve_{tt}^{\lambda})\big)\,\ve^{\lambda}_t\,\di x\,\di\tau&\leq&
\int_0^t\|T^{\lambda}\|^2_{L^2(\Omega;\S)}\,\di\tau +\|T^{\lambda}(t)\|^2_{L^2(\Omega;\S)}\nn\\[1ex]
&+&\|\sigma^{\lambda}(0)+\D(\ve_{t}^{\lambda}(0))\|^2_{L^2(\Omega;\S)}
+\int_0^t\|\ve_t^{\lambda}\|^2_{L^2(\Omega;\S)}\,\di\tau\nn\\[1ex]
&+&\nu\|\ve_t^{\lambda}(t)\|^2_{L^2(\Omega;\S)}+\tilde{C}(T)\,,
\end{eqnarray}
where the constant $\tilde{C}(T)$ does not depend on $\lambda$ and $\nu>0$ is any positive constant. Theorem \ref{tw:4.1} and (\ref{eq:4.13}) imply that the first, second, third and fourth term on the right hand-side of (\ref{eq:4.15}) are bounded. Moreover, the assumption $\nabla M(\dev(T(0)))\in L^2(\Omega;\SS)$ and the convexity of the function $M$ yield that
\begin{eqnarray}
\label{eq:4.16}
\int_{\Omega}M_{\lambda}(\dev(T(0)))\,\di x\leq \int_{\Omega}M(\dev(T(0)))\,\di x
\leq \int_{\Omega}\nabla M(\dev(T(0)))\,\dev(T(0))\,\di x\,.
\end{eqnarray}
Combining the inequalities (\ref{eq:4.13})-(\ref{eq:4.16}) the equality (\ref{eq:4.11}) becomes
\begin{eqnarray}
\label{eq:4.17}
&&\int_{\Omega}M_{\lambda}(\dev(T^{\lambda}(t)))\,\di x + \int_0^t\int_{\Omega}\D^{-1}T^{\lambda}_tT^{\lambda}_t\,\di x\,\di\tau  +\frac{1}{2}\int_{\Omega}\D(\ve^{\lambda}_{t}(t))\ve^{\lambda}_t(t)\,\di x\nn\\[1ex] 
&\leq&\tilde{C}(T)+\nu\|\ve_t^{\lambda}(t)\|^2_{L^2(\Omega;\S)}
\end{eqnarray}
where the constant $\tilde{C}(T)$ does not depend on $\lambda$. Choosing $\nu>0$ sufficiently small we complete the proof.$\mbox{}$ \hfill $\Box$\\[1ex]
{\bf Remark:}\hspace{2ex} The sequences $\{\theta^{\lambda}\}_{\lambda>0}$ is bounded in the space $L^2(0,T;H^1(\Omega;\R))$ and the sequence $\{\theta^{\lambda}_t\}_{\lambda>0}$ is bounded in the space $L^2(0,T;L^2(\Omega;\R))$, hence it contains a subsequence (again denoted using the superscript $\lambda$) such that $\theta^{\lambda}\rightarrow \theta$ a.e. in $\Omega\times(0,T)$. The continuity of $f$ and $T_{\frac{1}{\epsilon}}$ yield that 
$$f\big(T_{\frac{1}{\epsilon}}(\theta^{\lambda}+\tilde{\theta})\big)- f\big(T_{\frac{1}{\epsilon}}(\theta+\tilde{\theta})\big)\rightarrow 0 \quad \mathrm{a.e.\, in}\quad \Omega\times(0,T)$$ and $\big|\,f\big(T_{\frac{1}{\epsilon}}(\theta^{\lambda}+\tilde{\theta})\big)- f\big(T_{\frac{1}{\epsilon}}(\theta+\tilde{\theta})\big)\,\big|$ is bounded independently of $\lambda$. From the dominated Lebesgue theorem we conclude that for all $q\geq 1$   
$$f\big(T_{\frac{1}{\epsilon}}(\theta^{\lambda}+\tilde{\theta})\big)- f\big(T_{\frac{1}{\epsilon}}(\theta+\tilde{\theta})\big)\rightarrow 0\quad \mathrm{in}\quad L^q(0,T;L^q(\Omega;\R))\,.$$\\[1ex]
Theorem \ref{tw:4.1} and \ref{tw:4.2} yield that the sequence of stresses $\{T^{\lambda}\}_{\lambda>0}$ is bounded in $H^1(0,T;L^2(\Omega;\S))$. 
However, this information is not enough to pass to the limit in the system (\ref{eq:3.1}) as $\lambda$ tends to zero. We need to improve the convergence of the sequence $\{T^{\lambda}\}_{\lambda>0}$.
\begin{tw}$\mathrm{(Strong\, convergence\, of\, stresses)}$\\
\noindent
\label{tw:4.3}
Let us assume that the given data satisfy all requirements of Theorem \ref{tw:1.2}. Then,
\begin{eqnarray*}
\int_{\Omega}\D^{-1}(T^{\lambda}-T^{\mu})(T^{\lambda}-T^{\mu})\,\di x
\longrightarrow 0
\end{eqnarray*}
for $\lambda$,$\mu\rightarrow 0^+$ uniformly on bounded time intervals. 
\end{tw}
{\bf Proof:}\hspace{2ex} Calculate the time derivative
\begin{eqnarray}
\label{eq:4.18}
\lefteqn{\frac{d}{dt}\Big(\frac{1}{2}\int_{\Omega}\D\big(\ve^{\lambda}-\ve^{\mu} -(\ve^{p,\lambda}-\ve^{p,\mu})\big) \big(\ve^{\lambda}-\ve^{\mu} -(\ve^{p,\lambda}-\ve^{p,\mu})\big)\,\di x\Big)=}\nn\\[1ex] 
&=& \int_{\Omega}\D\big(\ve^{\lambda}-\ve^{\mu} -(\ve^{p,\lambda}-\ve^{p,\mu})\big) \big(\ve^{\lambda}_t-\ve^{\mu}_t -(\ve^{p,\lambda}_t-\ve^{p,\mu}_t)\big)\,\di x =\nn\\[1ex]
&=&\int_{\Omega}\D\big(\ve^{\lambda}-\ve^{\mu} -(\ve^{p,\lambda}-\ve^{p,\mu})\big) \big(\ve^{\lambda}_t-\ve^{\mu}_t\big)\,\di x \nn\\[1ex]
&-&\int_{\Omega}\Big(G^{\lambda}(\dev(T^{\lambda}))-G^{\mu}(\dev(T^{\mu}))\Big)
\Big(\dev(T^{\lambda})-\dev(T^{\nu})\Big)\,\di x\nn\\[1ex]
&=&\int_{\Omega}\big(\sigma^{\lambda}-\sigma^{\mu}+(\D(\ve^{\lambda}_t-\ve^{\mu}_t)(\ve^{\lambda}_t-\ve^{\mu}_t))\big)(\ve^{\lambda}_t-\ve^{\mu}_t)\,\di x - \int_{\Omega}\D(\ve^{\lambda}_t-\ve^{\mu}_t)(\ve^{\lambda}_t-\ve^{\mu}_t)\,\di x\nn\\[1ex]
&-&\int_{\Omega} \big(f\big(T_{\frac{1}{\epsilon}}(\theta^{\lambda}+\tilde{\theta})\big)-
f\big(T_{\frac{1}{\epsilon}}(\theta^{\mu}+\tilde{\theta})\big)\big) \big(\mathrm{div}\,u^{\lambda}_t-\mathrm{div}\,u^{\mu}_t\big)\,\di x \nn\\[1ex]
&-&\int_{\Omega}\Big(G^{\lambda}(\dev(T^{\lambda}))-G^{\mu}(\dev(T^{\mu}))\Big)
\Big(\dev(T^{\lambda})-\dev(T^{\nu})\Big)\,\di x\,.
\end{eqnarray}
Using the fact that the given data for two approximation steps are equal and integrating with respect to time, we conclude that
\begin{eqnarray}
\label{eq:4.19}
\lefteqn{\frac{1}{2}\int_{\Omega}\D\big(\ve^{\lambda}-\ve^{\mu} -(\ve^{p,\lambda}-\ve^{p,\mu})\big) (\ve^{\lambda}-\ve^{\mu} -(\ve^{p,\lambda}-\ve^{p,\mu}))\,\di x}\nn\\[1ex] 
&+&\int_0^t\int_{\Omega}\D\big(\ve^{\lambda}_t-\ve^{\mu}_t\big) (\ve^{\lambda}_t-\ve^{\mu}_t)\,\di x\,\di\tau \leq D\int_0^t\int_{\Omega} |f\big(T_{\frac{1}{\epsilon}}(\theta^{\lambda}+\tilde{\theta})\big)-
f\big(T_{\frac{1}{\epsilon}}(\theta^{\mu}+\tilde{\theta})\big)|^2\,\di x\,\di\tau\nn\\[1ex] 
&-&\int_0^t\int_{\Omega}\Big(G^{\lambda}(\dev(T^{\lambda}))-G^{\mu}(\dev(T^{\mu}))\Big)
(\dev(T^{\lambda})-\dev(T^{\nu}))\,\di x\,\di\tau\,,
\end{eqnarray}
where $D$ does not depend on $\lambda,\,\nu$. The remark before Theorem \ref{tw:4.3} and standard methods for maximal monotone operators finish the proof (cf. \cite{26}, \cite{23} and \cite{24}).$\mbox{}$ \hfill $\Box$\\[4ex]
{\large\em\bf Proof of Theorem \ref{tw:1.2}:}\hspace{2ex} From the definition of the Yosida approximation we obtain
\begin{eqnarray}
\label{eq:4.20}
\lefteqn{
\int_0^t\|J_{\lambda}(\dev(T^{\lambda}))-\dev(T)\|^2_{L^2(\Omega;\SS)}\di\tau\leq}\nn\\[1ex] &&\int_0^t\|J_{\lambda}(\dev(T^{\lambda}))-\dev(T^{\lambda})\|^2_{L^2(\Omega;\SS)}\di\tau
+\int_0^t\|\dev(T^{\lambda})-\dev(T)\|^2_{L^2(\Omega;\SS)}\di\tau\nn\\[1ex]
&=&\lambda\int_0^t\|\nabla M_{\lambda}(\dev(T^{\lambda}))\|^2_{L^2(\Omega;\SS)}\di\tau
+\int_0^t\|\dev(T^{\lambda})-\dev(T)\|^2_{L^2(\Omega;\SS)}\di\tau\,.
\end{eqnarray}
Theorem \ref{tw:4.2} implies that the sequence $\{\nabla M_{\lambda}(\dev(T^{\lambda}))\}_{\lambda>0}$ is bounded in $L^2(0,T;L^2(\Omega;\SS))$ and from Theorem \ref{tw:4.3} we conclude that the sequence $\{T^{\lambda}\}_{\lambda>0}$ is a Cauchy sequence in the space $L^{\infty}(0,T;L^2(\Omega;\S))$. Passing to the limit in (\ref{eq:4.20}) with $\lambda\rightarrow 0^+$ we have that the sequence $J_{\lambda}(\dev(T^{\lambda}))\rightarrow\dev(T)$ in $L^2(0,T;L^2(\Omega;\SS))$, hence it contains a subsequence such that $J_{\lambda}(\dev(T^{\lambda}(x,t)))\rightarrow\dev(T(x,t))$ for almost all $(x,t)\in (0,T)\times\Omega$. We know that 
$$\ve^{p,\lambda}_t(x,t)=\nabla M \big(J_{\lambda}(\dev(T^{\lambda}(x,t)))\big)= |J_{\lambda}(\dev(T^{\lambda}(x,t)))|^{p-1}J_{\lambda}(\dev(T^{\lambda}(x,t)))\,,$$
therefore
\begin{eqnarray}
\label{eq:4.21}
\ve^{p,\lambda}_t(x,t)\rightarrow |\dev(T(x,t))|^{p-1}\dev(T(x,t))\quad\mathrm{for\,\,a.\,a.}\quad (x,t)\in(0,T)\times\Omega\,.
\end{eqnarray}
Now, multiplying equation $(\ref{eq:3.1})_5$ by $v\in H^1(\Omega;\R)$ and integrating with respect to $\Omega$ we obtain
\begin{eqnarray}
\int_{\Omega}\theta_t^{\lambda}\,v\,\di x+\int_{\Omega}\nabla\theta^{\lambda}\nabla v\,\di x+ \int_{\Omega}f\big(T_{\frac{1}{\epsilon}}(\theta^{\lambda}+\tilde{\theta})\big) \mathrm{div}\,u_t^{\lambda}\,v\,\di x &=& \int_{\Omega}T_{\frac{1}{\epsilon}}(\dev(T^{\lambda})\,\ve^{p,\lambda}_t)\,v\,\di x\,.\nn
\end{eqnarray}
The information (\ref{eq:4.21}) implies  that
$$T_{\frac{1}{\epsilon}}\big(\dev(T^{\lambda}(x,t))\,\ve^{p,\lambda}_t(x,t)\big)- T_{\frac{1}{\epsilon}}\big(|\dev(T(x,t))|^{p+1}\big)\rightarrow 0\,\,\, \mathrm{for\,\,a.\,a.}\,\,\, (x,t)\in(0,T)\times\Omega\,$$
and moreover 
$$\big|T_{\frac{1}{\epsilon}}\big(\dev(T^{\lambda}(x,t))\,\ve^{p,\lambda}_t(x,t)\big)- T_{\frac{1}{\epsilon}}\big(|\dev(T(x,t))|^{p+1}\big)\big|\leq\frac{2}{\epsilon}\,.$$
From the Dominated Lebesgue theorem we conclude that
$$T_{\frac{1}{\epsilon}}\big(\dev(T^{\lambda})\,\ve^{p,\lambda}_t(x,t)\big)- T_{\frac{1}{\epsilon}}\big(|\dev(T)|^{p+1}\big)\rightarrow 0\,\,\, \mathrm{in}\,\,\, L^2(0,T;L^2(\Omega;\SS))\,.$$ 
Let us recall that (remark after Theorem \ref{tw:4.2})
$$f\big(T_{\frac{1}{\epsilon}}(\theta^{\lambda}+\tilde{\theta})\big)- f\big(T_{\frac{1}{\epsilon}}(\theta+\tilde{\theta})\big)\rightarrow 0\quad \mathrm{in}\quad L^{q}(0,T;L^{q}(\Omega;\R))$$
for all $q\geq 1$, hence the function
$$f\big(T_{\frac{1}{\epsilon}}(\theta^{\lambda}+\tilde{\theta})\big)v\rightarrow  f\big(T_{\frac{1}{\epsilon}}(\theta+\tilde{\theta})\big)v\quad \mathrm{in}\quad L^{2}(0,T;L^{2}(\Omega;\R))\,,$$
where $v\in H(\Omega;\R)$. Finally, the energy estimates (Theorem \ref{tw:4.1} $+$ Theorem \ref{tw:4.2}) give us the following information: The sequence $\{u^{\lambda},\ve^{p,\lambda}\}$ is bounded in $H^1(0,T;H^1(\Omega;\R^3)\times L^2(\Omega;\SS))$. The above informations are enough to pass to the limit in the Yosida approximation and get the solution in the sense of Definition \ref{de:1.1}. $\mbox{}$ \hfill $\Box$\\[1ex]
{\bf Remark: }\hspace{1ex} The regularity of $u_t\in L^2(0,T;H^1(\Omega;\R^3))$ implies that the function \\$\theta_t\in L^2(0,T;L^2(\Omega;\R))$. From the regularity theory for the parabolic equations we conclude that the heat equation in $(\ref{eq:1.1})_5$ is satisfied for almost all $(x,t)\in\Omega\times(0,T)$. Hence the solution from Definition \ref{de:1.1} is the $L^2$- strong solution - for the definition we refer to \cite{26}.\\[3ex]
{\bf Acknowledgments: }
This work was supported by the Grant of the Polish National Science Center 
OPUS 4 2012/07/B/ST1/03306.

\end{document}